\documentclass[12pt]{amsart}

\def\AA{\mathcal{C}}

\usepackage[square,compress,comma, numbers,sort]{natbib}
\usepackage[colorlinks=true, citecolor=red, linkcolor=blue]{hyperref}
\usepackage{amsfonts,mathtools}

\allowdisplaybreaks[4]

\usepackage{amssymb}
\usepackage{color}
\newcommand{\abs}[1]{\left\lvert #1 \right\rvert}

\DeclarePairedDelimiterXPP\pk[1]{\mathbb{P}}\{ \}{}{ #1}
\DeclarePairedDelimiterXPP\E[1]{\mathbb{E}}\{ \}{}{	#1}

\DeclarePairedDelimiterXPP\ind[1]{\mathbb{I}}( ){}{	#1}

\def\FRE{\mbox{Fr\'{e}chet }}

 \usepackage{xparse}

\NewDocumentCommand{\ceil}{s O{} m}{%
  \IfBooleanTF{#1} 
    {\left\lceil#3\right\rceil} 
    {#2\lceil#3#2\rceil} 
}
\NewDocumentCommand{\floor}{s O{} m}{%
  \IfBooleanTF{#1} 
    {\left\lfloor#3\right\rfloor}
    {#2\lfloor#3#2\rfloor}
}

\newcommand{\norm}[1]{\lVert #1 \rVert}
\def\x{\vk{x}}

\definecolor{c20}{rgb}{0.,0.7,0.}
\definecolor{c30}{rgb}{0.,0.,1.}
\definecolor{c40}{rgb}{1,0.1,0.7}
\definecolor{c50}{rgb}{1,0,0}
\definecolor{c60}{rgb}{1,0.9,0.1}
\definecolor{c70}{rgb}{0.50,1.00,0.00}

\def\cE#1{{\textcolor{c30}{#1}}}
\def\cE#1{#1}

\numberwithin{equation}{section}
\newtheorem{theo}{Theorem}[section]
\newtheorem{sat}[theo]{Proposition}
\newtheorem{de}[theo]{Definition}
\newtheorem{lem}[theo]{Lemma}

\newtheorem{example}[theo]{Example}
\newtheorem{korr}[theo]{Corollary}
\newtheorem{remark}[theo]{Remark}

\numberwithin{equation}{section}

\newcommand{\proofprop}[1]{\textsc{Proof of Proposition} \ref{#1}}
\newcommand{\prooflem}[1]{\textsc{Proof of Lemma} \ref{#1}}
\newcommand{\proofkorr}[1]{\textsc{Proof of Corollary} \ref{#1}}

\newcommand{\QED}{\hfill $\Box$}

\newcommand{\COM}[1]{}

\def\IF{\infty}

\newcommand{\R}{\mathbb{R}}
\newcommand{\inr}{\in \R}

\newcommand{\BQN}{\begin{eqnarray}}
\newcommand{\EQN}{\end{eqnarray}}
\newcommand{\BQNY}{\begin{eqnarray*}}
\newcommand{\EQNY}{\end{eqnarray*}}
\def\ldot{, \ldots,}

\def\todis{\overset{d}\rightarrow}

\newcommand{\kb}[1]{\boldsymbol{#1}}
\newcommand{\vk}[1]{\kb{#1}}

\def\bqny#1{\begin{eqnarray*} #1 \end{eqnarray*}}
\def\bqn#1{\begin{eqnarray} #1 \end{eqnarray}}

\newcommand{\BS}{\begin{sat}}
\newcommand{\ES}{\end{sat}}
\newcommand{\BT}{\begin{theo}}
\newcommand{\ET}{\end{theo}}
\newcommand{\BK}{\begin{korr}}
\newcommand{\EK}{\end{korr}}

\newcommand{\BEX}{\begin{example}}
\newcommand{\EEX}{\end{example}}

\newcommand{\BD}{\begin{de}}
\newcommand{\ED}{\end{de}}

\newcommand{\BIT}{\begin{itemize}}
\newcommand{\EIT}{\end{itemize}}
\newcommand{\BDI}{\begin{description}}
\newcommand{\EDI}{\end{description}}

\newcommand{\BRM}{\begin{remark}}
\newcommand{\ERM}{\end{remark}}

\newcommand{\BEL}{\begin{lem}}
\newcommand{\EEL}{\end{lem}}

\newcommand{\nelem}[1]{{Lemma \ref{#1}}}
\newcommand{\neprop}[1]{{Proposition \ref{#1}}}

\def\X{\vk{X}}

\newcommand{\equaldis}{\stackrel{{d}}{=}}
\newcommand{\eqZ}{\stackrel{{z.}}{=}}
\newcommand{\eqMZ}{\stackrel{{max-z.}}{=}}

\def\X{\vk{X}}

\def\TT{\mathcal{T}}

\def\tPr{\vk \Theta^{[h]} }

\def\TTT{\mathcal{T}}
\def\TT{\mathcal{T}}
\def\intT{\int_{T_0}}
\def\intT{\int_{\TT}}
\def\intTO{\sum_{s \in T_0} }

\begin{document}

\title{Multivariate Max-Stable Processes and Homogeneous Functionals}

\author{Enkelejd  Hashorva}
\address{Enkelejd Hashorva, Department of Actuarial Science \\
	Faculty of Business and Economics\\
University of Lausanne,\\
UNIL-Dorigny, 1015 Lausanne, Switzerland}
\email{Enkelejd.Hashorva@unil.ch}

\author{Alfred Kume }
\address{Alfred Kume,   Department of Statistics\\
	University of Kent, UK}
\email{alfred.kume@kent.ac.uk}

\bigskip

\date{\today}
 \maketitle

\begin{quote}

{\bf Abstract:} 
Multivariate max-stable processes are important for both theoretical investigations and various statistical applications motivated by the fact that these are limiting processes, for instance of stationary multivariate regularly varying time series, \cite{klem}.  In this contribution we  explore the relation between homogeneous functionals and multivariate max-stable processes and discuss  the connections between  multivariate max-stable process and  zonoid /  max-zonoid equivalence. We  illustrate our results considering Brown-Resnick and Smith  processes. 

\end{quote}
{\bf Key Words:} Multivariate max-stable process; homogenoeus functions; stationary process;   zonoid equivalence; max-zonoid equivalence;  Brown-Resncik process;  \\
{\bf AMS Classification:} Primary 60G15; secondary 60G70\\

\section{Introduction}  
Let $\vk X(t)=(X_1(t) \ldot X_d(t)), t\in \TT=\R^p$ be a $d$-dimensional  max-stable process with continuous  sample paths  and \FRE marginal distribution functions (df's); here $d,p$ are positive integers.  
In the light of  de Haan characterisation,  see e.g., \cite{deHaan,dom2016} 
we shall consider for simplicity $\vk X$ such that for some $\alpha>0$ it has the following representation (in distribution)
\bqn{\label{eq1}
	\vk X(t) =  \max_{i\ge 1} \Gamma_i^{-1/\alpha} \vk Z^{(i)}(t), \quad t\in \TT,
}
where $\Gamma_i= \sum_{k=1}^i E_k$ with $E_k, k\ge 1$ unit exponential random variables (rv's) independent of $\vk Z^{(i)}$'s, which are independent copies of a $d$-dimensional process $\vk Z(t)=(Z_1(t) \ldot Z_d(t)),t\in \TT$ with continuous  sample paths and non-negative components.
As in \cite{MolchanovSPA}  (therein $\vk Z_i$'s have  strictly positive components) see also \cite{MR3335107}, for any  $t_i \in \TT, 
\vk x_i\in (0,\IF)^d, i\le n$  
	\bqn{\label{eqfin}
(\pk{ \vk X(t_i) \le \vk x_i, 1 \le i\le n})^c&=& 
\pk{ \vk X(t_i) \le \vk x_i/c^{1/\alpha}, 1 \le i\le n} \notag\\
& =& \exp \Bigl( - c \E{ \max_{ 1 \le i\le d, 1 \le j \le n   }  Z_i^\alpha(t_j)/x_{ij}^\alpha    }  \Bigr)
}
is valid for any $c>0$. Eq. \eqref{eqfin}  is the so-called  max-stability property of $\vk X$; here we write often $\vk X$ instead of  $\vk X(t),t\in \TT$ and similarly for other processes and refer to $\vk Z$ as the spectral process of $\X$.
 
Write next $C $ for  the space of all  continuous functions  $ f: \TT \to [0,\IF) ^d $. It is well-known that $C$ can be equipped with a metric that turns it into a Polish space for which its  Borel $\sigma$-field coincides with $\AA$, the $\sigma$-field generated by projection maps $\pi_t, t\in T_0$, with $T_0$ a dense countable subset of $\TTT$. 

 Given $h\in \TT$ and some norm  $\norm{\cdot}$  on  $\R^d$,   define the tilted process $\vk \Theta^{[h]} $ (its law depends on the chosen norm $\norm{\cdot}$) by 
\bqn{ \label{fjalen}
	\pk{\vk  \Theta^{[h]} \in A } = \frac{1}{\E{ \norm{\vk Z(h)}^\alpha}} \E{ \norm{\vk Z(h)}^\alpha \mathbb{I}( \vk Z/\norm{\vk Z(h)} \in A) }
}	  for all  $A \in \AA$ with $\mathbb{I}(\cdot)$ the indicator function.   For notational simplicity we have assumed that $\vk \Theta^{hj]}$ is defined in the same probability space as $\vk Z$.\\
 In view of \eqref{eqfin},  if $R$ is a non-negative  rv  with $\E{R^\alpha}=1$ and independent of $\vk Z$, then clearly from \eqref{eqfin}   
\bqn{\label{RS}
	\widetilde{\vk Z}(t)=R \vk Z(t) 
} 
is another spectral process for $\vk X$;  in our notation $ a\vk x = (a x_1 \ldot a x_d), a\inr, \x \inr^d$. \\
An interesting alternative to random scaling is tilting, if for some $h\in \TT$ we have that $\pk{
	\max_{1 \le i \le n}  Z_i(h)>0}=1$.
Indeed, \eqref{eqfin} can be re-written as  (set $a_h:= \E{\norm{ \vk Z(h)}^\alpha}$ which is positive and finite by the assumption on \FRE marginals of $\vk X$)
\bqny{  
	\ln \pk{ \vk X(t_i) \le \vk x_i, 1 \le i\le n}
	& =& -  \E{\norm{\vk Z(h)}^\alpha/a_h \max_{ 1 \le i\le d, 1 \le j \le n   }  a_h( Z_i(t_j)/\norm{\vk Z(h)})^\alpha  /  x_{ij}^\alpha    } 
}
and thus  $\widetilde{\vk Z}(t)= a_h^{1/\alpha}\vk \Theta^{[h]}(t),t\in \TT$  is also another spectral process for $\vk X$. 
By definition of  $\vk  \Theta^{[h]}$  we have that 
\bqn{\label{lemBA}  
	\norm{\widetilde{\vk Z}(h)}^\alpha= a_h\norm{ \vk \Theta^{[h]}(h)}^\alpha =a_h \in (0,\IF)
} almost surely.  When $d=1$ in view of Balkema's Lemma given in \cite{HaanPickands}[Lem 4.1] 
(see Section 4 below for details) we have that any spectral process $ \widetilde{\vk Z}$ of $\vk X$ that 
satisfies \eqref{lemBA} has the same law as $a_h^{1/\alpha}\vk \Theta^{[h]}$. If 
$\pk{	\max_{1 \le i \le d}  Z_i(h)=0} \in (0,1)$, then 
$a_h^{1/\alpha}\vk \Theta^{[h]}$ can not be a spectral process for $\vk X$.  It is nonetheless possible to construct 
a spectral process for $\vk X$ by utilising a family of $\vk  \Theta^{[h]}$'s, see Section 2. 

\cE{It follows from \eqref{eqfin} that for given  two spectral processes $\vk Z$ and $\widetilde{\vk Z}$ of $\vk X$ 
	and all  maps $H( f) = 
	\max_{1 \le i \le d, 1 \le j \le n} ( f_i(t_j))^\alpha/x_{ij}, f\in C$ with  $t_j$'s in $\TT$ and $x_{ij}$'s in 
	$(0,\IF)$   
} 
\bqn{ \label{zr} 
	\E{  H(\vk Z)} = \E{ H(\widetilde{\vk Z})}. 
}
All maps  $H$ defined above  belong to the class $E_\alpha$  of  all  non-negative measurable $\alpha$-homogeneous maps $H:C\mapsto [0,\IF]$; here   $H$ is  $\beta$-homogeneous  means that
$H(c  f) = c^\beta H( f)$ for any $c > 0,  f\in C$.

In  \neprop{th0} below  we show how to construct $\tilde { \vk Z}$ and prove further that \eqref{zr}  holds for all spectral processes $\vk Z, \tilde{\vk Z}$ and all  $H\in E_\alpha$.  
It is known from \cite{MolchanovBE}  that homogeneous functions play a crucial role for the study of max-stable  random vectors. Therefore  the claimed validity of \eqref{zr} does not come as a surprise.

In Section 2 we  discuss briefly the implications of \eqref{zr} for stationary max-stable processes, whereas in  Section 3 we focus on the relations between zonoid equivalence, max-zonoid equivalence and homogeneous functions.  Section 4 is dedicated to Smith and Brown-Resnick  max-stable processes where 
 we derive also tractable formulas for their fidi's complementing previous results in \cite{MR3335107}. All the proofs are relegated to Section 5. 

To this end, we mention that numerous  results for max-stable processes and their representations exist in the literature, see for instance \cite{kab2009a,stoev2010max,WangStoev,MR3449778}. Our findings  in this paper, which have certain consequences for stationary max-stable processes,  are motivated by recent contributions  \cite{klem,BP,WS,HBernulli,bookSoulier} concerned with   multivariate regularly varying time series and their relation to  max-stable processes.

 \section{Results}
 
Let $\vk X$  be a max-stable process as in the Introduction with paths in $C$,  de Haan representation \eqref{eq1} and spectral process $\vk Z$  with sample paths in $C$ such that \eqref{14} holds.  
In view of \cite{dom2016}  for any compact set $K\subset \TT$ we have 
\bqn{ \label{wehave} 
	\E*{ \sup_{t\in K} \norm{\vk Z(t)}^\alpha}< \IF,
}	
which together with \cite{HBernulli}[Lem 7.1] implies that we can assume without loss of generality that 
\bqn{ \pk*{\sup_{t\in \TT} \norm{\vk Z(t)}>0}=1. 
	\label{14}
}
 For notational simplicity we shall suppose hereafter  that 
\bqn{ \label{imino} 
	\E{ \norm{\vk Z(t)}^\alpha}=1, \quad \forall t\in \TT.
}
In the following we consider $\vk \Theta^{[h]}$'s  to be independent and defined in the same non-atomic probability space  $(\Omega, \mathcal{F}, \mathbb{P})$, this is possible in view of \cite{Kallenberg}[Corr. 5.8]. Since $(C,\mathcal{C})$ is a Polish space and $\vk X$ has almost surely sample paths in $C$, in view of   \cite{MR100291}[Lem p.\ 1276] $\vk X$ can be realised also as a random process defined on $(\Omega, \mathcal{F}, \mathbb{P})$, which we shall assume in the sequel.

 Hereafter $T_0=\mathbb{Q}^p $  with $\mathbb{Q}$ the set of rational numbers,  $W$ is a $T_0$-valued   rv   with  probability mass 
 function (pmf) $p(t)> 0, t\in T_0$ being independent of everything else. 
  
\def\Nn{W}

We state next our first result: 

\BS \label{th0} Let $\vk X$ be max-stable with sample paths in $C$ and spectral process $\vk Z$ such that \eqref{14} and \eqref{imino} hold. Define  $\vk \Theta^{[h]}$ as in \eqref{fjalen} and let $\widetilde{\vk Z}$ be a random process with sample paths in $C$. \\ 
i) If $\widetilde{\vk Z}$ is a   spectral process of $\vk X$ such that $\pk{\sup_{t\in \TT} \norm{\widetilde{\vk Z}(t)}>0}=1$, then \eqref{zr} 
holds for all  $H\in E_\alpha$.\\ 
Conversely, if $\E{\norm{\vk Z(h)}^\alpha}=\E{\norm{\widetilde{\vk Z}(h)}^\alpha}=1,h\in \TT$ for some $\alpha>0$  satisfying \eqref{zr} with   $H=H_h( f)= \norm{ f(h)}^\alpha \Gamma(f), f\in \AA $  for all   $\Gamma\in E_0,h \in \TT$, then 
the corresponding max-stable processes $\vk X$ and $\widetilde{\vk X}$  of $\vk Z$ and $\widetilde{ \vk Z}$, respectively  are equal in law.  \\
ii) All random processes $ \vk Z_\Nn$ given by  
\bqn{ 
	\vk Z_\Nn(t)= \frac{1}{ (\intTO\norm{\vk \Theta^{[\Nn]}(s)}^\alpha p(s)  )^{1/\alpha} } \vk \Theta^{[\Nn]}(t), \quad t\in \TT
	\label{24}
}
are spectral processes for $\vk X$. 
\ES 

\BRM  
 i) In view of  \neprop{th0}  the law of $\vk \Theta^{[h]}, h\in \TT$ does not depend 
 on the particular choice of the spectral process $\vk Z$ since by \eqref{zr}
$$ \E{ \mathbb{I}( \vk \Theta^{[h]} \in A)} = \E{ \norm{\vk Z(h)}^\alpha  \mathbb{I}( \vk Z/\norm{\vk Z(h) \in A} )}= 
\E{ \norm{\widetilde{\vk Z}(h)}^\alpha  \mathbb{I}(  \widetilde{\vk Z}/\norm{\widetilde{\vk Z}(h)}  \in A )}, \quad \forall A \in \AA
$$
if  $ \widetilde{\vk Z}$ is another spectral process of $\vk X$ (recall we assume \eqref{imino}).  \\
ii) If $\vk X$ is max-stable stationary with sample paths in $C$, unit \FRE marginals and   spectral process $\vk Z$,  then by definition $B^h \vk X$ and $\vk X$ have the same law for any $h\in \TT$; here $B^h \vk X(t)= \vk X(t- h), h,t \in \TT$. In the light of  \eqref{eqfin} 
this is equivalent with $B^h \vk Z$ is a   spectral process for $\vk X$ for any $h\in \TT$, which in view of    \neprop{th0} implies for any $\widetilde{ \vk Z} $  a    spectral process of $\vk X$  
\bqn{  \label{t2} 
  \E{\Gamma_\alpha( \vk Z)}&=& 	\quad \E{\Gamma_\alpha( B^h\vk Z)} = \quad \E{\Gamma_\alpha( B^h \widetilde{\vk Z})}, \quad  \forall \Gamma_\alpha\in E_\alpha. 
} 
By \neprop{th0} the above also implies the stationarity of $\vk X$. Note in passing that our claim here  extends \cite{Htilt}[Thm 4.3,6.9] to the vector-valued setup.

\ERM

{\bf Example 1}: Let $\vk Z(t),t\in \TT$ be a BRs process as in \neprop{prop1}. 
For a given functional $F\in E_0$ which is shift-invariant in the sense that $F(B^h f) = F(f)$ for any $h\in \TT, f\in C$ define a new spectral process $\widetilde{\vk Z}(t)= \vk Z(t) F(\vk Z), t\in \TT$. Suppose that  $\E{\norm{\widetilde{\vk Z}(t_0)}^\alpha}\in (0,\IF)$ for some $t_0 \in \TT$. Since for given  $h\in \TT$ by \eqref{t2} 
\bqny{ 
	\E{  \Gamma_\alpha(\widetilde{\vk Z})}&=& 
	 \E{ \Gamma_\alpha ( \vk Z F (\vk Z))}= \E{ \Gamma_\alpha ( B^h\vk Z F (B^h \vk Z))} 
	 = \E{ \Gamma_\alpha( B^h \widetilde{\vk Z})} 
}
for all $ \Gamma \in E_\alpha$ we have that  $\widetilde{\vk Z}$ is also a BRs process.

 We set below $\vk \Theta:=\vk \Theta^{[0]}$ and present   three  equivalent conditions for the stationarity of a max-stable process $\vk X$ with sample paths in $C$ extending thus \cite{Htilt}[Thm 4.3].

\BK If  $\vk X, \vk Z$ are as in \neprop{th0}, then  $\vk X$  is stationary  if and only if:\\
 i) For any  $T_0$-valued rv  $\Nn $ with   pmf   $ p(t)>0,t\in T_0$ being independent of everything else 
 \bqn{\label{spectr}
 	\vk  Z_\Nn (t) = \frac{  1} { ( \intTO \norm{\vk \Theta(s-\Nn)}^\alpha   p(s) )^{1/\alpha}}  \vk \Theta (t-\Nn), \quad t\in \TT  
 }
 is a  spectral process for $\vk X$.\\
 ii) For all  $h \in \TT, \Gamma \in E_0$ and all positive integers  $k\le d$
\bqn{  \label{t2B}\E{   Z^\alpha _k(h) \Gamma(\vk Z)} = \E{ {  Z^\alpha_k(0)}  \Gamma(B^h \vk Z) 
		 }.
} 
iii) For all  $h \in \TT$ and all positive integers $k\le d$ 
\bqn{ \vk \Theta^{[h,k]}  \equaldis B^h \vk \Theta^{[0,k]},
\label{muha2}
} 
where $\vk \Theta^{[h,k]}  = (\Theta_{1}^{[h,k]}  \ldot \Theta_{d}^{[h,k]})$ is defined by 
$$\pk{\vk \Theta^{[h,k]} \in A}= 
\E{ \frac{ Z_k^\alpha (h)}{ \E{Z_k^\alpha (h)} } \mathbb{I}( \vk Z/ Z_k(h)   \in A)}, \quad  \forall A\in \AA,
$$
with  $ \vk \Theta^{[h,k]}(t),t\in \TT$   equal to $(1 \ldot 1) \in \R^d$  if $\E{Z_k^\alpha (h)}=0$.
\label{prop1}
\EK

\BRM \label{obvious} 
  \label{statio}  
i)  A general condition for the stationarity of  $\vk X$ in terms of spectral processes   based on the  findings of \cite{MR847381} is derived in \cite{MR2453345}[Thm 1] for the discrete setup.  When   $\vk Z(t), t\in \TT$ has strictly positive components for any $t\inr$,  a simple condition for the stationarity of $\vk X$ is given  in \cite{MolchanovSPA}.\\
ii) Using  \eqref{t2} it follows that  $\vk X$ is stationary if and only if 
\bqn{ \vk \Theta^{[h]} \equaldis B^h \vk \Theta^{[0]},
	\label{vogel}
}  
which is  first shown for $d=1$ in \cite{Htilt}[Thm 4.3, Eq. (4.6)].\\
iii) Stationary of max-stable processes can  be alternativly studied by relating it to the  shift-invariance of the corresponding tail measure   as in \cite{klem,PH2020}. \\
iv) If $\TT=\R$ we can similarly consider max-stable process with cadlag sample paths. All the results derived above remain unchanged, see \cite{PH2020} which deals with stationary $\vk X$. 
\ERM

\section{Zonoid and max-zonoid equivalence} In this section we discuss the relations between zonoid, max-zonoid equivalence and homogeneous functions. Essentially, we state some known results presenting some short and new proofs. \\ 
We shall consider below the case $\alpha=1$. Let $Z= (Z_1 \ldot  Z_n)$ and $Z^*=(Z_1^* \ldot  Z_n^*)$ be two random vectors not identically equal to zero with integrable components. As in \cite{MolchanovBE} we shall call
$Z$ and $Z^*$     zonoid equivalent if  
\bqn{\label{allu}
	\E{ F_{u}(Z)} = \E{ F_u(Z^*) }, \quad \forall u=(u_1 \ldot u_n) \inr^n,
}
where $F_u(Z)=\abs{\sum_{i=1}^n u_i Z_i } $;  abbreviate the above as $Z \eqZ Z^*$.\\
The distribution of $Z$, in view of Hardin \cite{Hardin}[Thm 1.1] is uniquely determined if we know  
$\E{ \abs{\sum_{i=0}^{n} u_i Z_i }}$ for any $u_i\inr, 0 \le i \le n$ assuming that  $Z_0=1$ almost surely. Consequently, if $(1, Z) \eqZ(1, Z^* )$, then $Z \equaldis Z^*$. 
\BEL If $Z$ and $Z^*$ have non-negative components, then $Z \eqZ Z^*$ if and only if for any $1$-homogeneous measurable function $H: \R^n \to [0,\IF) $ 
\bqn{\label{homE}
	\E{ H(Z)} = \E{H(Z^*)}.
}
\label{lemnosol}
\EEL 
As in \cite{MolchanovBE} we shall call non-negative $Z$ and $Z^*$ max-zonoid equivalent, if \eqref{allu} holds with $F_u(Z)= \max_{0 \le i \le n} u_i Z_i, $  where $u_0Z_0=0$. We use the abbreviation $Z\eqMZ Z^*$;  note in passing that $Z \eqMZ Z^*$  is the same as $X\equaldis X^*$, where $X$ and $X^*$ are max-stable random vectors whose spectral process are $Z$ and $Z^*$, respectively.  
The counterpart of Hardin's result is   Balkema's Lemma \cite{HaanPickands}[Lem 4.1], which implies  that $\E{ {\max_{0 \le i \le n} u_i Z_i }} $ for any $u_i\in [0,\IF) 
, 0 \le i \le n$  uniquely identifies the distribution of $(1,Z)$. Consequently, by Balkema's Lemma  
$$(1,Z) \eqMZ (1, Z^*)$$
is equivalent with the equality in distribution $Z \equaldis Z^*.$

\BEL If $Z$ and $Z^*$ have non-negative components, then $Z \eqMZ Z^*$ if and only if for any $1$-homogeneous measurable function $H: \R^n \to [0,\IF) $   \eqref{homE} holds. 
\label{lemnosolB}
\EEL

\BRM \label{raus} 
i) The claim of \nelem{lemnosol} holds for  $Z$ and $Z^*$ that can have negative components  requiring further that $H$ is an even function $\R^n \mapsto [0,\IF)$, which  follows from 
\cite{StoevWangSPL}[Thm 1.1] and  \cite{MolchanovBE}[Thm 2].  \\
ii)  A direct implication of \nelem{lemnosol} and \nelem{lemnosolB} is that 
$Z\eqMZ Z^* $ is equivalent with  $Z \eqZ Z^*$;  
see \cite{MolchDiag}[Thm 2.1], \cite{StoevWangSPL}[Thm 1.1]. 
\ERM

\section{Smith and Brown-Resnick processes}
In this section we consider briefly   Smith  and Brown-Resncik processes, see e.g., \cite{MolchanovSPA,MR3335107} for details. As shown in the aforementioned article these models are natural limiting models and therefore can be utilised for various statistical applications, which rely often on various tractable formulas for the fidi's of those processes. 
\subsection{Smith processes}
For a given parameter $\alpha>0$, we consider multivariate Smith processes that are   constructed by a given deterministic $[0,\IF)^d$-valued  function  $\vk L(t),t\in \TT$ with continuous  components $L_i,i\le d$     satisfying 
\bqn{ \label{dush} 
	 0<\intT \norm{\vk L(t) }^\alpha_* \lambda(dt)  \le \intT \sup_{s\in K} \norm{\vk L(t-s) }_*^\alpha \lambda(dt)    <  \IF
}
for any compact set $K \subset \R^p$ and some norm $\norm{ \cdot}_*$ on $\R^d$.
We define a multivariate max-stable Smith process $\vk X$ with  paths in $C$ by specifying a spectral process  $\vk Z$ of $\vk X$ as follows 
\bqn{ \label{mm} \vk Z(t) =(1/ p(\Nn))^{1/\alpha} \vk L(t- \Nn), \quad t\in \TT,
}
with $\Nn$ a $\TT$-valued  rv  with positive pdf  $ p(t)>0, t\in \TT$. For any norm $h\in \TT$ and any $\norm{\cdot}$ on $\R^d$ we have  that 
$	\vk \Theta^{[h]} $ calulcated with respect to this norm is given by 
\bqn{ \label{thash} 
	\vk \Theta^{[h]} (t)= \vk L(t-   h+S)/ \norm{\vk L( S)} = B^ h \vk \Theta^{[0]}(t), \quad t\in \TT,
}	  
with $S$ a  $\TT$-valued  rv   having  pdf $\norm{\vk L(t)} ^\alpha/c>0, t\in \TT$ with $c= \intT \norm{\vk L(t)} ^\alpha \lambda(dt)$.  \\
Note in passing that the inequality in \eqref{dush} is necessary in view of \eqref{wehave}. 
Since \eqref{thash} implies \eqref{vogel},  then the Smith max-stable process  $\vk X$ with paths in $C$ is stationary.

Set below    $c_k=\intT  {L_k^\alpha (t)}\lambda(dt) \in [0,\IF)$ for any positive integer $k\le d$. If $c_k>0$,  then   
\bqn{ \label{thash2} 
	\vk \Theta^{[h,k]}(t) =\vk L(t - h+ S_k)/ { L_k(S_k)}, \quad t\in \TT,
}	 
where  the $\TT$-vaéued  rv    $S_k $ has pdf ${ L_k^\alpha(t)}/c_k,t\in \TT$. 
Since both $\vk \Theta^{[h]}$ and $	\vk \Theta^{[h,k]}$ are known explicitly, we can calculate the fidi's  of $\vk X$  as shown next.  Hereafter $\norm{\vk x}_\IF= \max_{1 \le i \le d} \abs{x_i}, \vk x=(x_1 \ldot x_d)\inr^d$.

\BEL   \label{shtun}
Both \eqref{thash} and \eqref{thash2} are valid and for any $ t_i \in \TT,\vk x_i\in [0,\IF)^d, 1 \le i\le  n$, if further 
$c_k\in  (0,\IF), k\le d$  and $c_\IF=\intT \norm{\vk L(t)}_\IF ^\alpha \lambda(dt)$ we have  
\bqn{ \label{infargAGent} 
\lefteqn{ 	- \ln \pk{ \vk X(t_i) \le \vk x_i, 1 \le i\le n}} \notag \\
&	=	& c_\IF  \sum_{1\le l \le n} \E{  \norm{ \vk L ( S)/\vk x_{l}}_\IF^\alpha / \norm{\vk L(S)}_\IF^\alpha 
		\mathbb{I}( infargmax_{  1 \le j \le n }   \norm{ \vk L            (t_j- t_l+S)/\vk x_{j}}_\IF =l  )  }\\
	&=& \label{infargAGent2} 
	 \sum_{ (k,l) \in J} \frac{c_k}{x_{kl}^\alpha} 
	\pk*{	infargmax_{  (i,j) \in J  }  \frac{ L_i(t_j-t_l+S_k)}{x_{ij} {L_i(S_k)} } =(k,l)    },
}
where $J=\{(i,j), 1\le i \le d, 1 \le j \le p\}$ and   the infargmax funcitonal in \eqref{infargAGent2} is defined (also in the sequel) taking  the maximum  with respect to the lexicographical order in $\mathbb{Z}^2$.
 
\EEL

 \COM{
 
 For any $t_î\in \TT,x_i\in [0,\IF)^d, i\le  n$ 
 \bqn{ \label{infargAGent} 
 	- \ln \pk{ \vk X(t_i) \le \vk x_i, 1 \le i\le n}
 	&=&	  \sum_{ (k,l) \in J} \frac{c_k}{x_{kl}^\alpha} \pk{
 		infargmax_{  (i,j) \in J  }  \frac{ (L_i(t_j-t_l+S_k))^\alpha}{(x_{ij} {L_i(S_k)}) ^\alpha} =(k,l)    } \\
 	&=&  \sum_{ (k,l) \in J} \frac{c_k}{x_{kl}^\alpha} \pk{
 		infargmax_{  (i,j) \in J  }  \frac{ (L(t_j-t_l+S_k))^\alpha}{(x_{ij} {L(S_k)}) ^\alpha} =(k,l)    } \\
 	&=&  \sum_{ (k,l) \in J} \frac{c_k}{x_{kl}^\alpha} \pk{
 		(i,j) \in J:     \frac{ (L(t_j-t_l+S_k))^\alpha}{(x_{ij} {L(S_k)}) ^\alpha} \le   
 		\frac{ (L(S_k))^\alpha}{(x_{kl} {L(S_k)}) ^\alpha}    } \\
 	&=&  \sum_{ (k,l) \in J} \frac{c_k}{x_{kl}^\alpha} \pk{
 		(i,j) \in J:     \frac{ (L(t_j-t_l+S_k))^\alpha}{(  {L(S_k)}) ^\alpha} \le   
 		\frac{  x_{ij}^\alpha}{x_{kl}  ^\alpha}    } \\
 	&=&  \sum_{ (k,l) \in J} \frac{c_k}{x_{kl}^\alpha} \pk{
 		(i,j) \in J:    \frac{\sigma_k^2}{\sigma_i^2} W^2 -  \frac{1}{\sigma_i^2 }(t_j-t_l+ \sigma_k W)^2 \le 2 \ln 
 		\frac{  x_{ij}}{x_{kl} }    } \\
 	&=&  \sum_{ (k,l) \in J} \frac{c_k}{x_{kl}^\alpha} \pk{
 		(i,j) \in J:    \sigma_k^2 W^2 -  (t_j-t_l+ \sigma_k W)^2 \le 2 \ln 
 		\frac{  x_{ij}}{x_{kl} }    } \\
 	&=&  \sum_{ (k,l) \in J} \frac{c_k}{x_{kl}^\alpha} \pk{
 		(i,j) \in J:    \sigma_k  (t_j-t_l) W  - (t_j- t_l)^2/2 \le  \ln 
 		\frac{  x_{ij}}{x_{kl} }    } \\
 }
}
 {\bf Example 2}: Assume that $\alpha=1$ and  $\vk L(t), t\in \TT$ has deterministic components  $ L_i, i\le d$ being equal to the pdf of a Gaussian random vector 
 $\vk Y=(Y_1 \ldot Y_p)$ with independent $N(0,1)$ components.  In view of \nelem{shtun}  for 
 $t_i=(t_{i1} \ldot t_{ip}) \in \TT,\vk x_i\in [0,\IF)^d, 1 \le i\le  n$ we have 
 \bqny{  
 	\lefteqn{	- \ln \pk{ \vk X(t_i) \le \vk x_i, 1 \le i\le n}} \notag \\
 	&=&	  \sum_{ (k,l) \in J} \frac{1}{x_{kl}} \pk*{   (i,j) \in J  \setminus \{(k,l)\}: \  (  \sum_{r=1 }^p   (t_{lr} -t_{jr}) Y_r - 
 		(t_{jr} - t_{lr})^2/2 )  \le  \ln (x_{ij} /x_{kl}) } ,
 }
where  $J= \{ (i,j), 1 \le j \le d, 1 \le j \le p\}.$

\subsection{Brown-Resnick processes}   
In the Brown-Resnick model,  $\alpha=1$ and the spectral process  is the log-Gaussian one  
\bqn{\label{grillI}
	\vk Z(t)= e^{ \vk Y(t) - \vk V(t) /2}, \quad \vk V(t)=(Var(Y_1(t) \ldot Var(Y_d(t)), \quad t\in \TT,
}
where $\vk Y(t),t\in \TT $ is a $\R^d$-valued Gaussian process with centered components. 
\BEL  \label{tru}  If $\vk X$ is a Brown-Resnick   max-stable process as defined in \eqref{grillI}, then the law of $\vk X$ depends only on the matrix-valued function   with $ij$th entry equal $\gamma_{ij}(t,s)=Var(Y_i(t)- Y_j(s)),s,t\in \TT$ where $1 \le i,j \le d$.   In particular, if $\gamma_{ij}(t,s), i,j\le d$ depend only on the difference $t-s$, then $\vk X$ is stationary.  
Further, for any  $t_i\in  \TT ,\vk x_i \in (0,\IF)^d, 1 \le i\le n$
\bqn{ \label{infargAG} 
	\lefteqn{	- \ln \pk{ \vk X(t_i) \le \vk x_i, 1 \le i\le n}} \notag \\
	\hspace{1 cm}&=&	  \sum_{ (k,l) \in J } \frac{1}{x_{kl}} \pk{ \forall  (i,j) \in J \setminus \{ (k,l) \}:  Z_i(t_j) - Z_k(t_l)   - Var(  Z_i(t_j) - Z_k(t_l) )/2    \le \ln (x_{ij}/x_{kl} )},
}
provided that $(\vk Y(t_1) \ldot \vk Y(t_n))$ possesses a density. 
\EEL

\BRM   
A similar  formula to \eqref{infargAG} is shown for the special case that $\vk X$ has stationary increments in \cite{MR3335107};  that  assumption is not needed in our case. It is worth noting that the formula \cite{Nik} 
derived for H\"usler-Reiss distributions appears in \cite{Nik}. Note further that  the other claims of \nelem{tru}  are stated without proof in \cite{MR2453345}.
\ERM

\section{Proofs}

\proofprop{th0} $i$) 
Let $\vk Z$ and $\widetilde{\vk Z}$ be two spectral processes of $\vk X$  with sample paths in $C$ and let $H_h( f)= \norm{ f(h)}^\alpha  \Gamma(f/\norm{f(h)}), f\in C, h\in \TT$ for some norm $\norm{\cdot}$ on $\R^d$ and $\Gamma$ a measurable functional $C\mapsto [0,\IF]$. Considering   the case that  $\norm{\vk x}=\norm{\vk x}_\IF=  \max_{1 \le i \le d} \abs{x_i} $ we have 
\bqn{\label{yyyshi}  
	\E{ H_h(\vk Z)} =\E{ \norm{ \vk Z(h)}^\alpha  \Gamma(\vk Z/ \norm{ \vk Z(h)})} =  \E{  \Gamma(\tPr) }, \quad 
	\forall h \in \R ,
}  
with $\tPr$ defined in \eqref{fjalen}.  Since $\Gamma$ can be approximated by simple functions,  the claim in \eqref{zr}  for $H=H_h$  follows by showing that the law of $\tPr$ is uniquely defined by that of $\vk X$ and this does not depend on the spectral process $\vk Z$. 
We shall  assume for simplicity that 
$$ \E{\norm{\vk Z(h)}^\alpha } =1, \quad \alpha =1.$$
By the definition    $\norm{\tPr(h)}=\max_{1 \le j \le d} \Theta^{[h]}_j(h)=1$ almost surely. 
Using \eqref{eqfin} for any $s>1, \vk x_1 \ldot \vk x_n \in (0,\IF)^d, t_i \in \TT, i\le n $ and for any $u>0$ we have 
\bqny{ \lefteqn{\pk{ \norm{\vk X(h)} \le   s u, \vk X(t_i) \le  u\vk x_i , 1 \le i\le n \lvert  \norm{\vk X(h)} >  u  }}\\
 	&=& \frac{ u\pk{ \norm{\vk X(h)} \le   su, \vk X(t_i) \le u \vk x_i,  1 \le  i\le n}- \pk{ \norm{\vk X(h)} \le  u , \vk X(t_i) \le  u\vk x_i , 
		1 \le 	i\le n}}{u 
		\pk{\norm{\vk X(h)} >    u  }}\\
	&\to & \E{ \max( \norm{\vk Z(h)},  \max_{1 \le j \le d, 1 \le i \le n } Z_{j}(t_i)/x_{ij} ) - \max(  \norm{\vk Z(h)} /s, \max_{1 \le j \le d, 1 \le i \le n } Z_{j}(t_i)/x_{ij})}\\
	&=& \E{ ( \norm{\vk Z(h)}- \max( \norm{\vk Z(h)} /s, \max_{1 \le j \le d, 1 \le i \le n } Z_{j}(t_i)/x_{ij})_+}\\
	&=& \E{ ( 1- \max(1 /s, \max_{1 \le j \le d, 1 \le i \le n }  \Theta^{[h]}_{j}(t_i)/x_{ij})_+}\\
	&=& \E{ ( 1- \max_{1 \le j \le d, 0 \le i \le n }  \Theta^{[h]}_{j}(t_i)/x_{ij})_+}\\
	&=&  \pk{ R \le s,   R\vk{\Theta}^{[h]}(t_i)\le \vk x_i, 0 \le i \le n		}
}
as $u\to \IF$, with  
 $R$ a Pareto  rv  with survival function $1/r, r>1$ independent of $\tPr$ and we set $x_{0j}=1/s, t_0=h$. Consequently, we have  
 the convergence in distribution
$$ \Bigl ( \frac{\vk X(t_1)}{  \norm{\vk X(h) }} \ldot \frac{\vk X(t_n)}{ \norm{\vk X(h)}} \Bigr) \Bigl \lvert ( \norm{\vk X(h)} > n)   \todis 
(\tPr(t_1) \ldot \tPr(t_n)), \quad n\to \IF.$$
Hence  the fidi's of $\tPr$ are uniquely determined by those of $\vk X$. Hence  \eqref{yyyshi} holds and further  for all $\Gamma \in E_0, h\in \TT$ 
\bqn{\label{women} \E{\norm{\vk Z(h)} _\IF^\alpha \Gamma (\vk Z) }=\E{\norm{ \widetilde{\vk Z}(h)} _\IF^\alpha \Gamma (\widetilde{\vk Z}) }.    
} 
Note in passing that the above convergence in distribution follows also from \cite{BP} for any norm on $\R^d$.\\
Next consider a general norm $\norm{\cdot}$ on $\R^d$.  Since $ \E{\norm{ \vk Z(h)}} \in (0,\IF)$ for any $h\in \TT$  implies by the equivalence of the norms on $\R^d$,  then also $ \E{\norm{ \vk Z(h)}_\IF} \in (0,\IF)$ 
and by \eqref{14}, \eqref{imino} we can find a non-negative measurable function $q(t),t\in \TT$ such that 
$$S(\vk Z)= \intT \norm{ \vk Z(t)}_\IF  q(t) \lambda(dt) \in (0,\IF) $$ almost surely, with $\lambda(dt)$ the Lebesgue measure on $\TT$. Moreover, by the assumption on $\widetilde{\vk Z}$ we also have $\pk{S(\widetilde{\vk Z}) \in (0,\IF)}=1$. Consequently,  
for any  $H \in E_\alpha$, i.e., $H: C \mapsto [0,\IF]$ is an $\alpha$-homogeneous  measurable functional Fubini-Tonelli theorem yields 
\bqny{ \E{ H(\vk Z)} &=& \E*{ H(\vk Z)  \frac{S(\vk Z)}{S(\vk Z)}} = 
	  \intT \E*{  \norm{\vk Z(t)}_\IF \frac{H(\vk Z)  }{S(\vk Z)}} q(t) \lambda(dt) \notag \\
	&=:&  \intT \E{  \norm{\vk Z(t)}_\IF   F(\vk Z) } q(t) \lambda(dt)  \notag  \\
	&=& 
	\intT \E{  \norm{\widetilde{\vk Z}(t)}_\IF   F(\widetilde{\vk Z}) } q(t)\lambda(dt) \notag =  \E{  H(\widetilde{\vk Z}) }, 
}  
where we used the fact that $F\in E_0$  and \eqref{women} for the derivation of the second last equality above.\\
We show below the converse, i.e.,  for a given norm  $\norm{\cdot}$ on $\R^d$ we assume that \eqref{zr} holds for any $H_h= \norm{f(h)}^\alpha  \Gamma(f), h \in \TT, \Gamma \in E_0$  and prove that the   max-stable processes $\vk X$ and $\widetilde{\vk X}$ with spectral processes  $\vk Z$ and $\widetilde{\vk Z}$, respectively have the same fidi's. For given $t_1 \ldot t_n \in \TT$ let  
$H(\vk Z) =  \max_{ 1 \le i\le n}  \norm{ \vk Z(t_i)/\vk x_{i}}_\IF,
$ with  $\vk x_1 \ldot \vk x_n\in (0,\IF)^d$.    Putting $S(\vk Z)= \sum_{i=1}^n  \norm{ \vk Z(t_i)}_\IF$ which is positive whenever $H(\vk Z)$ is positive, then as above 
\bqny{ \E{ H(\vk Z)} = \E*{ H(\vk Z)  \frac{S(\vk Z)}{S(\vk Z)}} = \sum_{i=1}^n \E*{  \norm{\vk Z(t_i)}_\IF \frac{H(\vk Z)  }{S(\vk Z)}}	= 
	\sum_{i=1}^n \E{  \norm{\widetilde{\vk Z}(t_i)}_\IF   F(\widetilde{\vk Z}) }=  \E{  H(\widetilde{\vk Z}) }
}  
since $F \in E_0$ by the assumption on $H$. Consequently,  in view of  \eqref{eqfin}  $\vk X$ and $\widetilde{\vk X}$ have the same fidi's.\\

$ii)$  It follows easily that (recall $\norm{\vk \Theta^{[h]}(h)}=1$ almost surely for all $h \in \TT$)
$$I(W)= \intTO\norm{\vk \Theta^{[W]}(s)}^\alpha p(s)  \in (0,\IF)$$ almost surely and thus $\vk Z_W$ is well-defined and with sample paths in $C$. 
 In order to establish the proof we utilise the claim of statement $i)$. Let therefore $H\in E_\alpha$ be given. Using again Fubini-Tonelli  
 theorem we have 
\bqny{ 
	\E{ H( \vk Z_\Nn)}
	&=& \E*{\sum_{x\in T_0}  \frac{1}{ I(x)  }
		  H( \vk \Theta^{[x]} )     p (x)  } \\
	&=& \sum_{x\in T_0 }\E*{		 \norm{\vk Z(x)}^\alpha  \frac{1}{\intTO (\norm{\vk Z (s) } /\norm{\vk Z (x) }) ^\alpha p(s)  }
  H( \vk Z / \norm{\vk Z (x) }  )     p(x)  }\\
	&=& \E*{   H( \vk Z)    
		\sum_{x\in T_0} \frac{\norm{\vk Z(x)  }^\alpha} {\intTO \norm{\vk Z (s) } ^\alpha p(s)}    p(x)} 	=\E{   H( \vk Z) } 
} 
establishing the proof.  \QED  

\proofkorr{prop1} 
 $i)$ Assume that $\vk \X$ is stationary, i.e., \eqref{t2} holds and recall that $\norm{\vk \Theta^{[0]}(0)}=\norm{\vk \Theta(0)}=1$ almost surely for all $h\in \TT$.   It follows   that   
 \bqny{ \label{jake} 
 	J(\vk \Theta)= \intTO \norm{  B^W \vk \Theta  (s)}^\alpha p (s) \in (0,\IF)
 }
 almost surely, where  $\Nn$ is a $T_0$-valued   rv  with pdf $ p(t)>0, t\in T_0$ being further independent of $\vk Z$. 
 Next, for $\vk  Z_\Nn$ defined in \eqref{spectr},   
 using  \eqref{t2}  for the derivation of the third last equality,
  for all $F\in E_0, h \in T_0$   we have
 \bqny{ \E{ \norm{\vk  Z_\Nn(h)}^\alpha F(\vk Z_{\Nn}) } 
 	&=& \E*{  \sum_{x\in T_0}  
 		\norm{   \vk Z(0)}^\alpha 	\norm{   B^x \vk Z(h)}^\alpha  \frac{ 1}
 		{J(\vk Z)} 		 F(     B^x \vk Z)   p (x)    }\\
 	&=:& \sum_{x\in T_0}  \E{    \norm{   \vk Z(0)}^\alpha \Gamma( B^x \vk Z) }  p (x)  \\
 	&=& \sum_{x\in T_0}  \E{    \norm{   \vk Z(x)}^\alpha  \Gamma(\vk  Z) }  p(x) \\
 	&=&  \E*{  \norm{   \vk Z(h)}^\alpha  F(\vk Z)   \sum_{x\in T_0}   \frac{ \norm{   \vk Z(x)}^\alpha    p (x) }{   \intTO 
 			\norm{   \vk Z(s)}^\alpha    p(s) } }  \\
 	&=&  \E{  \norm{   \vk Z(h)}^\alpha  F(\vk Z)}.
 } 
Since $T_0$ is dense in $\TT$  and $\vk Z$ has sample paths in $C$ and thus it is stochastic continuous,  then by the above equality and \eqref{wehave},  applying the dominated convergence theorem  yields  for all
 $h\in \TT, k>0 $
$$ \E{ \norm{\vk  Z_\Nn(h)}^\alpha F(\vk Z_{\Nn}) \ind{F(\vk Z_{\Nn})\le k } }  =
\E{  \norm{   \vk Z(h)}^\alpha  F(\vk Z) \ind{F(\vk Z_{\Nn} )\le k }}.
$$ 		
Letting $k$ to infinity, by   \neprop{th0} we conclude  that $\vk  Z_\Nn$ is a spectral process for $\vk X$. Note in passing that if  $F \in E_0$, then  also $ F(f) \ind{F( f)\le k },f\in C$ belongs to $E_0.$\\
 Conversely, if $ \vk  Z_\Nn$ 
   defined in \eqref{spectr} is a spectral process for $\vk X$ with  $W$ independent of $\vk \Theta $ with pdf $ p(t)>0, t\in T_0$, then 
 for any $h \in T_0$  
 \bqny{ B^h \vk  Z_\Nn &=&
 	B^{\Nn+h}\vk \Theta  (t) \frac{  1}{ (  \intTO \norm{ B^\Nn \vk \Theta (s)}^\alpha  p(s) )^{1/\alpha}} \\
 	&\equaldis & B^{\Nn_h}\vk \Theta(t) 
 	\frac{ 1 }{ ( \intTO \norm{ B^{W_h} \vk \Theta  (s)}^\alpha p_h (s)  )^{1/\alpha}} 	=: \vk Z_{\Nn_h}, 
 }
 where $\Nn_h=\Nn+h$ has pdf  $p_h(t)=  p(t-h)$. Hence   $B^h \vk  Z_\Nn$ is a spectral process for $\vk X$  since it is equal in law with $\vk Z_{\Nn_h},$ which by the assumption is a spectral process for $\vk X$. Using the stochastic continuity of $\vk Z_W$ and \eqref{wehave}, then 
 $\vk Z_{\Nn_h}$ is a spectral process of  $\vk X$ for all $h \in \TT$,  thus $\vk X$ is stationary establishing  the proof. \\
  If \eqref{muha2} holds, then 
 by the definition \eqref{t2B} follows. If the latter is satisfied, then  for $\norm{\vk x}_\alpha = (\sum_{1 \le i\le d} \abs{ x_i}^\alpha )^{1/\alpha}$ and all $\Gamma \in E_0$ 
\bqn{\label{dalim2}
	 \E{ \norm{ \vk Z(t+h)}^\alpha_\alpha \Gamma(\vk Z) } = \E{ \norm{ \vk Z(t)}^\alpha_\alpha \Gamma(B^h\vk Z) }, \quad \forall h,t\in \TT .
	}
Since $\norm{\vk x}_\alpha$ vanishes if and only if $\norm{\vk x}$ vanishes, where $\norm{\cdot}$ is some norm on $\R^d$,  as in the proof of \neprop{th0} for given $t_i  h \in \TT, \vk x_n\in (0,\IF)^d, 1 \le i\le n$ and   
$$ H(\vk Z) =  \max_{ 1 \le i\le n}  \norm{ \vk Z(t_i+h)/\vk x_{i}}_\IF
$$ 
 by \eqref{dalim2} we have with 
$S_\alpha(\vk Z)= \sum_{i=1}^n  \norm{ \vk Z(t_i)}_\alpha $, which is positive whenever $H(\vk Z)$ is positive
\bqny{ \E{ H(\vk Z)} &=& \sum_{i=1}^n \E*{  \norm{\vk Z(t_i+h)}_\alpha^\alpha  \frac{H(\vk Z)  }{S_\alpha(B^{-h}\vk Z)}}
			=: 	\sum_{i=1}^n \E{  \norm{{\vk Z}(t_i+h)}_\alpha^\alpha    F(  \vk Z)}\\
			&=&
	\sum_{i=1}^n \E{  \norm{{\vk Z}(t_i)}_\alpha^\alpha    F( B^h \vk Z)} =  \E{  H( B^h \vk Z) }
}  
since $F \in E_0$ by the assumption on $H$. Consequently,   $\vk X$ and $B^h \X$ have the same fidi's and thus $\vk X$ is stationary. 
\QED

 \prooflem{lemnosol} First note that $Z \eqZ Z^*$ is equivalent with $$ (S, Z) \eqZ (S^*, Z^*),$$ 
 where $S= \sum_{i=1}^n Z_i$ and $S^*=\sum_{i=1}^n Z_i^*$ and $\E{S}= \E{S^*}= a$.  The assumption that  $Z$ and $Z^*$ are non-negative,  not identical to zero almost surely  and with integrable components  yields $a \in (0,\IF)$. Assume for simplicity that $a=1$. For  any  $H: \R^n \mapsto  [0,\IF)$ a  1-homogeneous, measurable function that  vanishes  at the origin 
 \bqn{\label{situa}
 	\E{H(Z)}= \E{H(Z ) \ind{S>0}}=\E{H(Z/S)  S} = \E{H(\widetilde Z)},
 }
 where $\widetilde Z$ is the tilted random vector with respect to $S$. From  $Z \eqZ Z^*$ we have that $\widetilde Z \eqZ \widetilde {Z^*}$ and consequently by the definition of zonoid equivalence 
 $$ ( \sum_{i=1}^n \widetilde Z_i, \widetilde Z) \eqZ  (\sum_{i=1}^n \widetilde{Z^*}_i, \widetilde{Z^*} ).$$
 Since by the definition $\sum_{i=1}^n \widetilde Z_i = \sum_{i=1}^n \widetilde{Z^*}_i=1$ almost surely, then 
 Hardin's result yields that $\widetilde Z \equaldis \widetilde{Z^*}$ and consequently using \eqref{situa}, by  the measurability of  $H$   
 \bqn{\label{kA} 
 	\E{H(Z)}= \E{H(\widetilde Z)}=\E{H(\widetilde{Z^*})}=\E{H(Z^*)},
 }	
 hence the proof is complete.  \QED
 
 \prooflem{lemnosolB}   
 First note that $Z \eqMZ Z^*$ is equivalent with 
 $ (M, Z) \eqZ (M^*, Z^*),$  where $M= \max_{1 \le i\le n} Z_i$ and $M^*=\max_{1 \le i \le n} Z_i^*$ and $\E{M}= \E{M^*}=a$. Since $\E{M} \le \sum_{i=1}^n \E{Z_i}$ and  $Z$ and $Z^*$ are non-negative integrable and not zero almost surely we have that $a\in (0,\IF)$. Suppose for simplicity that $a=1$. For  any  $H: \R^n \mapsto [0,\IF)$ a  1-homogeneous, measurable function that  vanishes  at the origin 
 $$\E{H(Z)}= \E{H(Z ) \ind{M>0}}=\E{H(Z/M)  M} = \E{H(\widetilde M)},$$
 where $\widetilde Z$ is the tilted random vector (tilted with respect to $M$). From the above and $Z \eqMZ Z^*$ we have that $\widetilde Z \eqMZ \widetilde {Z^*}$ and consequently 
 $$ ( \max_{1 \le i \le n} \widetilde Z_i, \widetilde Z) \eqMZ  (\max_{1 \le i \le n} \widetilde{Z^*}_i, \widetilde{Z^*} ).$$
 Since by the definition $\max_{1 \le i \le n} \widetilde Z_i = \max_{1 \le i \le n} \widetilde{Z^*}_i=1$ almost surely, then 
 Balkema's Lemma  yields that $\widetilde Z \equaldis \widetilde{Z^*}$ and consequently, the measurability of  $H$ implies \eqref{kA}, 
 hence the proof is complete.  \QED

 \prooflem{shtun} We give first the proofs of  \eqref{thash} and \eqref{thash2}. For any $h\in T$ by the definition of $\vk Z$  in 
 \eqref{mm}, i.e.,  $ \vk Z(t)= (1/ p (\Nn))^{1/\alpha} B^\Nn \vk L(t),  t\in \TT$, 
 with $W$ a $\TT$-valued  rv  with positive pdf  $ p(t)>0, t\in \TT$ using further the translation invariance of the Lebesgues measure on $\TT$ for all $h \in \TT$ 
 $$ \E{ \norm{\vk Z(h)}^\alpha }= \intT \norm{\vk L(t)}^\alpha \lambda(dt) =c.$$
 Since  $\intT \norm{\vk L(t)}^\alpha_* \lambda(dt)\in (0,\IF)$, by the equivalence of the norms on $\R^d$ we have that $c\in (0,\IF)$. 
  Given  $ A\in \AA$ 
 for any $h\in \TT$  Fubini-Tonelli theorem implies  
 \bqny{ \pk{\vk{ \Theta}^{ [h] } \in A} 
 	&=& \frac{1}{c}\E{ \norm{\vk Z(h)}^\alpha  \mathbb{I} (\vk Z/\norm{\vk Z(h)} \in A ) }\\
 	&=&  \frac{1}{c}\intT  \norm{(B^{-t} \vk L)(h)}^\alpha  \mathbb{I} (B^{-t} \vk L/\norm{(B^{-t} \vk L)(h)} \in A) \lambda(dt)\\
 	&  = & \intT  \mathbb{I}( B^{h-s}\vk L /\norm{ \vk L(s)}  \in A)  \norm{\vk L(s)}^\alpha/c  \lambda(ds)\\
 	& = &\pk{ B^{h-S}\vk L /\norm{\vk L(S)} \in A}, 
 }
 where the $d$-dimensional random vector   $ S$ has pdf  $\norm{\vk L(t)}^\alpha/c,t\in \TT$, 
 hence \eqref{thash} follows.

Since for any positive integer $k\le d$ by the translation invariance of the Lebesgue measure on $\TT$    
 $$ \E{ { Z_k^\alpha(h)} }= \intT { L_k^\alpha(t)}  \lambda(dt)=c_k $$
 for all  $A \in \AA, h\in \TT$,  if further $c_k>0$    
 \bqny{ \pk{\vk \Theta^{[h,k]} \in A}&=& \E{ { Z_k^\alpha(h)} \mathbb{I}(\vk Z/{ Z_k(h)} \in A) /c_k}\\
 	&  = & \intT  \mathbb{I}( B^{-t}\vk L /   {L_k(t)}  \in A)  L_k^\alpha(t) /c_k  \lambda(dt)\\
 	&= &\pk{ B^{h-S_k}\vk L/ { L_k(S_k)} \in A}, 
 }
 where the $\TT$-valued  random vector  $ S_k$  has  pdf  $ { L_k^\alpha (t)}/c_k,t\in \TT$, hence  
 \eqref{thash2} follows.\\
Next we show that  the fidi's of a  max-stable process $\vk X$ as in the Introduction can be determined  by  infargmax functional in terms of $\vk \Theta^{[k]}$, where the maximum and minimum are taken with respect to the lexicographical order.   In view of \eqref{eqfin} for any $t_i\in \TT, \vk x_i \in (0,\IF)^d, 1 \le i\le n$ we have (set  $\norm{\x}=: \norm{\vk x}_\IF=\max_{1\le i \le d} \abs{x_i}$)
\bqny{ 
	\lefteqn{	- \ln \pk{ \vk X(t_i) \le \vk x_i, 1 \le i\le n}} \notag \\
	&=&  \E{  \max_{1 \le j \le n }  \norm{ \vk Z(t_j)/\vk x_{j}}^\alpha    }\notag  \\
	&=&	  \sum_{1\le l \le n} \E{   \norm{ \vk Z(t_l)/\vk x_{l}}^\alpha \
		\mathbb{I}( infargmax_{  1 \le j \le n }   \norm{ \vk Z(t_j)/\vk x_{j}}^\alpha =l  )  } \notag \\
	&=&	  \sum_{1\le l \le n}\E{\norm{ \vk Z(t_l)}^\alpha } \E{   \frac{ \norm{ \vk Z(t_l)}^\alpha }{\E{\norm{ \vk Z(t_l)}^\alpha} } 
		\frac{ \norm{ \vk Z(t_l)/\vk x_{l}}^\alpha}{\norm{ \vk Z(t_l)}^\alpha } 
		\mathbb{I}( infargmax_{  1 \le j \le n }   \norm{ \vk Z(t_j)/\vk x_{j}} =l  )  } \notag \\
	&=&	c  \sum_{1\le l \le n} \E{  \norm{ \vk \Theta^{[t_l]}(t_l)/\vk x_{l}}^\alpha 
		\mathbb{I}( infargmax_{  1 \le j \le n }   \norm{ \vk \Theta^{[t_l]}           (t_j)/\vk x_{j}} =l  )  },
}
hence \eqref{infargAGent} follows from \eqref{thash}. Next,  if    $c_{k}=\E{Z_k^\alpha(t)}, t\in \TT$ 
is positive for any $k\le d $, setting $J=\{ (i,j), 1 \le j \le d, 1 \le j \le p\}$ for all $t_i \in \TT, \vk x_i \in (0,\IF)^d, 1 \le  i\le n$ we can write 
 using  \eqref{eqfin} 
 \bqn{ \label{franga}
 	\lefteqn{	- \ln \pk{ \vk X(t_i) \le \vk x_i, 1 \le i\le n}} \notag  \\
 	&=&  \E{ \max_{ (i,j) \in J}  Z_i^\alpha(t_j)/x_{ij}^\alpha    }\notag  \\
 	&=&	  \sum_{  (k,l) \in J} \frac{1}{x_{kl}^\alpha} \E*{ Z_{k}^\alpha(t_l) 
 		\mathbb{I} \Bigl( infargmax_{  (i,j) \in J  } \bigl( Z_i^\alpha(t_j)/x_{ij}^\alpha \bigr) =(k,l) \Bigr )  } \notag \\
 	&=&	  \sum_{  (k,l) \in J} \frac{\E{Z_{k}^\alpha(t_{l})}}{x_{kl}^\alpha} \E*{ \frac{Z_{k}^\alpha(t_l)}{\E*{Z_{k}^\alpha(t_l)}}
 		\mathbb{I} \Bigl( infargmax_ {  (i,j) \in J  } \bigl( Z_i^\alpha(t_j)/x_{ij}^\alpha\bigr) =(k,l)  \Bigr)  } \notag \\
 	&=&	  \sum_{( k,l) \in J} \frac{\E{Z_{k}^\alpha(t_{l})}}{x_{kl}^\alpha} \pk{
 		infargmax_{  (i,j) \in J  }  ( \Theta_i^{[t_l,k]}(t_j)/x_{ij})  =(k,l)    }.
 }
Note in passing that the above calculations  hold also when some $c_k$'s are equal to zero.  Applying the above formula and utilising  further \eqref{thash2} establishes the proof.  
 \QED

 \prooflem{tru} 
 In view of \cite{Htilt}[Lem 6.1] we have that $\vk \Theta^{[t_l,k]} $ is again log-Gaussian, but there is an additional deterministic trend function which can be calculated for each component separately as therein.  Specifically,  for any positive integer  $k\le d$ and all $h\in \TT$  we have
 \bqn{\label{jas}
\quad \quad  \quad ( \Theta^{[h,k]}_1(t) \ldot \Theta^{[h,k]}_d(t))   \equaldis \Big(  e^{  Y_1(t)- Y_k(h)    -
	 \frac{Var( Y_1(t)- Y_k(h) )}2 }\ldot
 	e^{  Y_d(t)- Y_k(h)    -\frac{Var( Y_d(t)- Y_k(h) )}2 }  \Bigr),
 }
 hence in view of \eqref{franga}  
 the fidi's  of $\vk X$ depends only on the  matrix-valued functions $P^{h,k}(t,s)$ with $ij$th entry equal to $p^{h,k}_{ij}(t,s)=cov(Y_i(t)- Y_k(h), Y_j(s)- Y_k(h)).$ Since
 $$p^{h,k}_{ij}(t,s) = [ \gamma_{ik}(t,h) + \gamma_{jk}(s,h)- \gamma_{ij}(t,s) ]/2, \quad 1 \le i,j\le d, h,s,t\in \TT,$$ 
 with  $\gamma_{ij}(t,s)=Var(Y_i(t)- Y_j(s))$, then the law of $\vk X$ depends only on $\gamma_{ij}$'s. If $\gamma_{ij}(t,s)$ depends only on the difference $t-s$, then 
 $$p^{h+a,k}_{ij}(t+a,s+a) = [ \gamma_{ik}(t,h) + \gamma_{jk}(s,h)- \gamma_{ij}(t,s) ]/2, \quad 1 \le i,j\le d, a,h,s,t\in \TT$$ 
and hence  by  \eqref{franga} we have that $(\vk X(t_1) \ldot \vk X(t_n))$ has the same law as $(\vk X(t_1+a) \ldot \vk X(t_n+a))$ for any $a\in \TT$ implying  $\vk X$ is stationary. Using \eqref{franga} the formula \eqref{infargAG} follows easily, hence the proof is complete. \QED

 { \bf Acknowledgement}:  We are in debt to the reviewers for their comments and suggestions which lead to improvements of the manuscript.  
 Partial 	support from SNSF Grant 200021-196888  is kindly acknowledged.

 \bibliographystyle{ieeetr}
 \bibliography{EEEA}
\end{document}